\documentclass[11pt]{amsart}
\usepackage{bm}
\usepackage{amssymb}
\usepackage{graphicx}	
\usepackage{xcolor}
\usepackage{amsmath}
\usepackage{enumerate}
\usepackage{blindtext}
\usepackage{multirow}
\usepackage{mathrsfs}

\newtheorem{thm}{Theorem}[section]
\newtheorem{Def}[thm]{Definition}
\newtheorem{rem}[thm]{Remark}
\newtheorem{cor}[thm]{Corollary}

\newtheorem{ex}[thm]{Example}
\newtheorem{prob}[thm]{Problem}
\newtheorem{lem}[thm]{Lemma}

\newcommand{\bdfn}{\begin{Def} \rm}
\newcommand{\edfn}{\end{Def}}
\newcommand{\tFAE}{The following are equivalent}
\newcommand{\tfae}{the following are equivalent}

\newcommand{\ra}{\rightarrow}
\newcommand{\Ra}{\Rightarrow}

\newcommand{\ci}{\subseteq}

\newcommand{\al}{\alpha}

\newcommand{\de}{\delta}
\newcommand{\e}{\varepsilon}

\newcommand{\si}{\sigma}
\newcommand{\Si}{\Sigma}
\newcommand{\ga}{\gamma}
\newcommand{\Ga}{\Gamma}

\newcommand{\mb}{\mathbb}
\newcommand{\mc}{\mathcal}

\newcommand{\iy}{\infty}
\newcommand{\om}{\omega}
\newcommand{\Om}{\Omega}

\newcommand{\beqa}{\begin{eqnarray*}}
\newcommand{\eeqa}{\end{eqnarray*}}
\newcommand{\vertiii}[1]{{\left\vert\kern-0.25ex\left\vert\kern-0.25ex\left\vert #1 
    \right\vert\kern-0.25ex\right\vert\kern-0.25ex\right\vert}}
\newcounter{cnt1}
\newcounter{cnt2}
\newcounter{cnt3}
\newcounter{cnt4}
\newcommand{\blr}{\begin{list}{$($\roman{cnt1}$)$} {\usecounter{cnt1}
\setlength{\topsep}{0pt} \setlength{\itemsep}{0pt}}}
\newcommand{\blR}{\begin{list}{\Roman{cnt4}.\ } {\usecounter{cnt4}
\setlength{\topsep}{0pt} \setlength{\itemsep}{0pt}}}
\newcommand{\bla}{\begin{list}{$(\alph{cnt2})$} {\usecounter{cnt2}
\setlength{\topsep}{0pt} \setlength{\itemsep}{0pt}}}
\newcommand{\bln}{\begin{list}{$($\arabic{cnt3}$)$} {\usecounter{cnt3}
\setlength{\topsep}{0pt} \setlength{\itemsep}{0pt}}}
\newcommand{\el}{\end{list}}

\begin{document}
\title[Unique Hahn-Banach extensions in Bochner $L^p$ spaces]{Unique Hahn-Banach extensions in Bochner $L^p$ spaces}
\author[Karak]{Sainik Karak}
\address[Sainik Karak]
{Department of Mathematics\\ I. I. T. Hyderabad\\ Telangana\\India,\textit{E-mail~:}\textit{ma22resch11001@iith.ac.in}}
\author[Paul]{Tanmoy Paul}
\address[Tanmoy Paul]
{Department of Mathematics\\ I. I. T. Hyderabad\\ Telangana\\India,\textit{E-mail~:}\textit{tanmoy@math.iith.ac.in}}

\subjclass[2000]{Primary 46A22, 46B20, 46E15 Secondary 46B22, 46M05 \hfill \textbf{\today} }

\begin{abstract}
Let $Y$ be a closed subspace of a real Banach space $X$. Let $(\Omega,\Si,\mu)$ be a complete finite measure space and let $1<p<\iy$. We prove that $Y$ has property-$(U)$, property-$(SU)$, property-$(wU)$ in $X$ if and only if $L^p(\mu,Y)$ has the corresponding property in $L^p(\mu,X)$. No Radon-Nikod\'{y}m assumption on $X^*$ is required. The proofs use $w^*$-measurable
representations of dual functionals and a simultaneous scalar-lifting construction. A norm-minimal representation for an element in $L^p(\mu,X)^*$ is found suitable in this case. We also study the point of continuity of the identity map under different topologies on the unit sphere of $L^p(\mu,X)^*$. For Bochner-integrable functionals, the pointwise $w^*-w$-continuity yields the uniqueness of the extension to $L^p(\mu,X^{**})$. Examples show why arbitrary $w^*$-measurable representatives cannot replace the norm-minimal representatives in this context. 
\end{abstract}

\maketitle

\section{Introduction}

\subsection{Preliminaries}
Uniqueness of Hahn–Banach extensions is closely related to the geometry of dual unit balls and the ideal structure of Banach spaces. This area is closely connected to the point of continuity of functionals on the dual unit sphere under different topologies, the asymptotic norming properties of the dual spaces, and the ideal structure of the underlying spaces (see \cite{DPR, DMR, SD, HL1, RS}). In \cite{HL} and, more recently, in \cite{DPR}, the authors have taken the initiative to investigate how well these properties are preserved under various sums and products of Banach spaces. In this study, refinements of the results under weaker assumptions are  presented. 
Throughout this paper, $X$ represents a real Banach space, and by a subspace we mean a closed subspace. By $NA(X)$ we mean the set of all functionals on $X$ which attain their norms on the unit ball of $X$. By $B_X$ and $S_X$, we denote the closed unit ball and the unit sphere of $X$, respectively.

\bdfn
Let $Y$ be a subspace of a Banach space $X$. Then,
\bla
\item\cite{Phelps}
 $Y$ is said to have {\it property-$(U)$} in $X$ if for each $f\in Y^*$, there exists unique $\widetilde{f}\in X^*$ such that $\|f\|=\|\widetilde{f}\|$ and $\widetilde{f}|_Y=f$.
 \item\cite{O} 
 $Y$ is said to have {\it property-$(wU)$} in $X$ if for each $f\in NA(Y)$, there exists unique $\widetilde{f}\in X^*$ such that $\|f\|=\|\widetilde{f}\|$ and $\widetilde{f}|_Y=f$.
 \item\cite{Rao} 
 $Y$ is said to be an {\it ideal} in $X$ if there exists a projection $P:X^*\to X^*$ such that $\ker(P)=Y^\perp$ and $\|P\|=1$.
\item\cite{O} 
$Y$ is said to have {\it property-$(SU)$} in $X$ if $Y$ has property-$(U)$ in $X$ and $Y$ is an ideal in $X$. 

\el
\edfn

Denote $Y^\#=\{f\in X^*:\|f\|=\|f|_Y\|\}$ for a subspace $Y$ of $X$. $Y^\#$ is a subset of $X^*$ but  not necessarily a subspace of $X^*$. The following characterizations of the notions mentioned above are useful in this study. Here we assume $Y$ and $X$ are as above.

\begin{thm}\cite{Lim, Phelps}
\tFAE.
\bla
\item $Y$ has property-$(U)$ in $X$.
\item $Y^\perp$ is Chebyshev in $X^*$.
\item For $y_1^*,y_2^*\in Y^\#$ satisfying $y_1^*+y_2^*\in Y^\perp$, then $y_1^*+y_2^*=0$.
\el
\end{thm}
 
\begin{thm}\cite{Lim, O}\label{T6}
\tFAE. 
\bla
\item $Y$ has property-$(SU)$ in $X$.
\item $Y^\#$ is a linear subspace of $X^*$.
\item For $y_1^*,y_2^*,y_3^*\in Y^\#$ satisfying $y_1^*+y_2^*+y_3^*\in Y^\perp$, we have $y_1^*+y_2^*+y_3^*=0$.
\el
\end{thm}

\begin{thm}\label{CwU}\cite{DP, Phelps}
\tFAE.
\bla
\item $Y$ has property-$(wU)$ in $X$. 
\item For $x^*\in X^*$, $P_{Y^\perp}(x^*)=\{0\}$ holds, whenever $x^*(y)=\|x^*\|$ for some $y\in S_Y$.
\item For $x_1^*,x_2^*\in NA(X)$, $x_1^*+x_2^*=0$ holds, whenever for $i=1, 2$, $x_i^*(y_i)=\|x_i^*\|$ for some $y_i\in S_Y$ and $x_1^*+x_2^*\in Y^\perp$.
\el 
\end{thm}

\bdfn\label{D1}\cite{HWW, HL}
A Banach space $X$ is said to be {\it Hahn-Banach smooth (weakly Hahn-Banach smooth)} if $X$ has property-$(U)$ (property-$(wU)$) in $X^{**}$ under the canonical embedding.
\edfn

The following is a well-known characterization of the phenomenon stated in Definition~\ref{D1}.

\begin{thm}\cite{HWW, HL1}\label{T7}
\tFAE.
\bla 
\item $X$ is Hahn-Banach smooth.
\item Each $x^*\in S_{X^*}$ has a unique norm preserving extension to a functional on $X^{**}$.
\item The identity mapping $I:(B_{X^*},w^*)\to (B_{X^*},w)$ is continuous on $S_{X^*}$.
\el
Moreover, an $x^*\in S_{X^*}$ has unique norm preserving extension to a functional on $X^{**}$ if and only if the identity mapping $I:(B_{X^*},w^*)\to (B_{X^*},w)$ is continuous at $x^*$.
\end{thm}

\subsection{Aim and objectives for this study}
Our first attempt is to answer the following in its full generalization.

\begin{prob}\label{Pr3}
Let $1<p<\iy$ and $(\Omega,\Sigma,\mu)$ be a probability space. Furthermore, let $Y$ be a subspace of a Banach space $X$ such that $Y$ has property-$(U)$ in $X$. Does $L^p(\mu,Y)$ have property-$(U)$ in $L^p(\mu,X)$?
\end{prob}

We recall the following from \cite{DPR} in this context. 
\begin{thm}\label{T9}
Let $(\Omega,\Sigma,\mu)$ be a probability space where $X^*$ has the RNP, then $Y$ has property-$(U)$ in $X$ if and only if $L^p(\mu,Y)$ has property-$(U)$ in $L^p(\mu,X)$ for $1<p<\iy$.
\end{thm}

In the absence of RNP, one can conclude that for those $f\in L^q(\mu,Y^*) (\subseteq L^p(\mu,Y)^*)$  has a Hahn-Banach extension $g\in L^q(\mu,X^*) (\subseteq L^p(\mu,X)^*)$ must be unique if $Y$ has property-$(U)$ in $X$.

It is clear that $Y$ has property-$(U)$ in $X$ if $L^p(\mu,Y)$ has property-$(U)$ in $L^p(\mu,X)$; we do not need the assumption of RNP for this case. If $X^*$ does not have the RNP, then we answer Problem~\ref{Pr3} in an affirmative sense.

A major concern of the spaces of type $L^p(\mu,X)$ is that we need to deal with the equivalence classes where two distinct representatives differ in a measure zero set. Different pointwise behaviour of the representatives concerning the properties stated above may cause difficulties in concluding about the equivalence class. A unified approach is adopted in this study to overcome such issues. This study may be useful in investigating other geometric properties of the spaces of Bochner-integrable functions.

Preliminaries and references on concepts related to Bochner integrable functions are included in the next section. 
A norm-preserving extension of a functional defined on its superspace is referred to as a {\it Hahn-Banach extension}. The limitations of this investigation are articulated in terms of the questions. Unless otherwise stated, when it comes $q$ for a given $p$ where, $1<p<\iy$ we mean $\frac{1}{p}+\frac{1}{q}=1$.

\subsection{Summary of the manuscript's architecture}
A brief review of the dual description of elements in $L^p(\mu,X)^*$ is given in section~2. It is demonstrated that for $\psi\in L^p(\mu,X)^*$ there exists a suitable $w^*$-measurable $\Phi:\Omega\ra X^*$ which satisfies $\psi(f)=\int_\Omega\langle f(\omega),\Phi(\omega)\rangle d\mu(\omega)$ for $f\in L^p(\mu,X)$.
In section~3 we discuss the aforesaid notions for a subspace of type $L^p(\mu,Y)$ in $L^p(\mu,X)$ under suitable assumptions on $Y$ in $X$. The above mentioned representations corresponding to $\varphi\in L^p(\mu,Y)^*$ and corresponding to its Hahn-Banach extension over $L^p(\mu,X)^*$ are found pointwise Hahn-Banach extension.  
In section~4, we discuss two examples related to the above mentioned representations. A necessary and sufficient condition for the Hahn-Banach smoothness of $L^p(\mu,X)$ is the equivalence of the topologies of $w^*$ and $\si(L^q(\mu,X)^*, L^p(\mu,X^{**}))$ on unit sphere of $L^q(\mu,X^*)$, is obtained in section~5. This extends the familiar characterisation of Hahn-Banach smoothness as discussed in Theorem~\ref{T7}.

\section{Applications of dual representations and scalar lifting maps}
Let $X$ be a Banach space, and $(\Omega, \Sigma, \mu)$ be a measure space. We follow the standard notation $L^p(\mu, X)$ to denote the space of all {\it $p$th-Bochner integrable functions}. These are $X$-valued strongly measurable functions defined on $\Omega$ satisfying $\int_\Omega \|f(\omega)\|^pd\mu(\om)<\iy$ for $1 \leq p< \infty$ and $\mbox{ess sup }\{\|f(\omega)\|:\omega\in \Omega\}<\iy$ for $p=\iy$. We recall that for a measure space $(\Omega,\Sigma,\mu)$ and a Banach space $X$, $g:\Omega\ra X$ is said to be {\it strongly measurable} if $g$ is a pointwise limit of a sequence of simple functions. Note that in this case the range of $g$ is necessarily separable.
Recall from \cite{BR, DU} that if $1/p + 1/q = 1$, $1 < p, q < \infty$, the space $L^q(\mu, X^*)$ is isometrically isomorphic to a norming subspace of $L^p(\mu, X)^*$ and that they coincide if and only if $X^*$ has the {\it Radon-Nikod\'{y}m property}. We use the abbreviation RNP to indicate the Radon-Nikod\'{y}m property. We refer the reader to the monographs \cite{CM, WL} for the necessary background on Bochner-integrable functions. \cite{DU} is a standard reference for the notion of RNP in Banach spaces.

We first recall the following theorems from \cite{CM} related to the duals of $L^p(\mu,X)$, for $1< p<\infty$. In the following, we assume that $(\Omega,\Sigma,\mu)$ is a finite measure space.

\begin{thm}\label{T10}\cite[Theorem 1.5.4]{CM}
Let $1\leq p<\iy$ and $q$ be the conjugate of $p$, then for each $\Ga\in L^p(\mu,X)^*$ there exists a function $\Psi:\Om\to X^*$ satisfying the following.
 \bla 
\item $\Psi$ is $w^*$-measurable.
\item The function $\omega\to \|\Psi(\om)\|\;($which we simply denote $\|\Psi(.)\|)$ is measurable and belongs to $L^q(\mu)$.
\item $\Ga(f)=\int_{\Omega}\langle f(\omega),\Psi(\omega)\rangle d\mu(\omega)$ for all $f\in L^p(\mu,X)$ and
\item $\|\Ga\|=\|\|\Psi(.)\|\|_q$
 \el
\end{thm}

For $1< p,q<\iy$, we introduce the space $\mathcal{L}_{w^*}^q(\mu,X^*)$, consisting of all $w^*$ -measurable functions $\Psi$ defined on underlying measure space for which there exists $g$, a scalar valued measurable function which satisfies $\|\Psi (\omega)\|\leq |g(\omega)|$ where $\int_\Omega |g(\omega)|^qd\mu(\omega)<\infty$. Given two functions $\Psi_1, \Psi_2\in \mathcal{L}_{w^*}^q(\mu,X^*)$, we say $\Psi_1\sim \Psi_2\iff$ for every $x\in X$, $\Psi_1(\omega)(x)=\Psi_2(\omega)(x)$ $\mu$-a.e. This induces an equivalent relation on $\mathcal{L}_{w^*}^q(\mu,X^*)$. The equivalence classes of the functions in $\mathcal{L}_{w^*}^q(\mu,X^*)$ represent the space $L_{w^*}^q(\mu,X^*)$. By $\mc{L}_\iy(\mu)$ we mean the vector space of all scalar-valued essentially bounded functions defined on $\Omega$.
We refer the reader to \cite[Subsection 1.5]{CM} for further details related to this notion. 

The identification stated in Theorem~\ref{T10} induces an isometric isomorphism 
$I:L_{w^*}^q(\mu,X^*)\to L^p(\mu,X)^*$ given by, 
\begin{equation}\label{Eq3}
I([\Psi])(f)=\int_\Om \langle f(\om),\Psi_0(\om)\rangle d\mu(\om)
\end{equation}
for $f\in L^p(\mu,X)$, $\Psi_0$ is any function in $[\Psi]$. 
The {\it quotient norm} of an arbitrary $[\Psi]\in L^q_{w^*}(\mu,X^*)$ is given here.

\[
\|[\Psi]\|_{q,w^*}=\inf\{\|a\|_q:\Phi\sim \Psi, a\in L^q(\mu), a\geq 0, \|\Phi(\om)\|\leq a(\om),\mu-a.e.\} 
\]

Recall the following result, which we also need in our subsequent discussion.

\begin{thm}\label{T1}
Let $(\Omega,\Sigma,\mu)$ be a finite measure space, and let $1< p<\iy$, then the bounded linear operator $I$ defined in equation~\ref{Eq3} is an isometric isomorphism onto $L^p(\mu,X)^*$, and for each class $[\Psi]\in L_{w^*}^q(\mu, X^*)$ there exists $\Psi_0\in [\Psi]$ such that the function $\omega\mapsto \|\Psi_0(\omega)\|$ (which we simply denote by $\|\Psi_0(.)\|$) is measurable and belongs to $L^q(\mu)$, and $\|[\Psi]\|_{q,w^*}=\|\|\Psi_0(.)\|\|_q$.
\end{thm}

We call such a representative stated in Theorem~\ref{T1} a {\it norm-minimal representative}. 

By {\it lifting} in $\mc{L}^\iy(\mu)$ we mean a linear map $\rho:\mc{L}^\iy(\mu)\to \mc{L}^\iy(\mu)$
that satisfies the conditions stated in Definition~1.5.1 in \cite{CM}. Following this definition, one can conclude that
$|\rho(f)|\leq \rho (|f|)$ everywhere for $f\in \mc{L}^\iy(\mu)$. To see this, first observe that $-|f|\leq f\leq |f|$. This leads to $\rho(-|f|)\leq \rho(f)\leq \rho(|f|)$ and this with the linearity of $\rho$ follows that  $|\rho(f)|\leq \rho (|f|)$.

If $[\Phi]\in L_{w^{*}}^q(\mu,X^{*})$ and $\Phi_1,\Phi_2$ are two norm-minimal representatives of $[\Phi]$, then
$$
\|\Phi_1(\om)\|=\|\Phi_2(\om)\|
\quad \mu\text{-a.e.}
$$
However, it need not follow that $\Phi_1(\om)=\Phi_2(\om)$ outside a single $\mu$-null set. Indeed, the equivalence relation only ensures that, for each fixed $x\in X$,
$$
\Phi_1(\om)(x)=\Phi_2(\om)(x)
\quad \mu\text{-a.e.},
$$
where the null set may depend on $x$. By contrast, if $[g]\in L^q(\mu)$ and $f_1,f_2\in[g]$, then $f_1(\om)=f_2(\om)$ $\mu$-a.e., and hence $|f_1(\om)|=|f_2(\om)|$, outside a single $\mu$-null set. In this work, our goal is to pass property-$(U)$ of $Y$ in $X$ to property-$(U)$ of $L^p(\mu,Y)$ in $L^p(\mu,X)$. The first natural difficulty that arises is as follows.

\begin{prob}\label{P1}
Let $[\Psi]\in L_{w^*}^q(\mu,Y^*)$ and $[\Phi]\in L_{w^*}^q(\mu,X^*)$ be a norm-preserving extension of $[\Psi]$. Also let $\Psi_0$ and $\Phi_0$ are the corresponding norm-minimal representatives of $[\Psi]$ and $[\Phi]$ respectively. Can we say $\Phi_0(\om)$ is a norm-preserving extension of $\Psi_0(\om)$ $\mu$-a.e.?
\end{prob}

\begin{rem}\label{R1}
Initially, we know that $\|\Phi_0(\om)\|=\|\Psi_0(\om)\|$ $\mu$-a.e. and for each fixed $y\in Y$, $\Phi_0(\om)(y)=\Psi_0(\om)(y)$ $\mu$-a.e. However, this is not sufficient because the null set changes as $y$ varies. To address this, we use a scalar lifting $\rho:\mc{L}_\iy(\mu)\to \mc{L}_\iy(\mu)$ to find norm-minimal representatives $\widetilde{\Phi_0}$ and $\widetilde{\Psi_0}$ with $[\widetilde{\Phi_0}]=[\Phi_0]$ and $[\widetilde{\Psi_0}]=[\Psi_0]$ satisfying $$\widetilde{\Phi_0}(\om)|_Y=\widetilde{\Psi_0}(\om)\quad \mu-a.e.$$  But for the new representatives, we do not know whether $\|\widetilde{\Phi_0}(w)\|=\|\widetilde{\Psi_0}(\om)\|$ $\mu$-a.e. The proof is complete if we can ensure that $\|\widetilde{\Phi_0}(\om)\|=\|\Phi_0(\om)\|$ $\mu$-a.e and $\|\widetilde{\Psi_0}(\om)\|=\|\Psi_0(\om)\|$ $\mu$-a.e., since we already have $\|\Phi_0(\om)\|=\|\Psi_0(\om)\|$ $\mu$-a.e. The following lemma is introduced to overcome this. However, one can see the details of this remark later in Theorem~\ref{TU}.
\end{rem}
  
For the following lemma, we fix $1<p<\iy$ and $q=\frac{p}{p-1}$.

\begin{lem} \label{LL}
Let $[\Phi]\in L_{w^*}^q(\mu, X^*)$, and choose a norm-minimal representative $\Phi_0$ of $[\Phi]$. 
Then there exists another norm-minimal representative $\widetilde{\Phi}$ of $[\Phi]$ such that $[\Phi_0]=[\widetilde{\Phi}]$ satisfies the following. \bla
\item $\widetilde{\Phi}(\om)\in X^*$, for all $\om$
\item $\|\Phi_0(\om)\|=\|\widetilde{\Phi}(\om)\|$, for $\mu$-a.e.
\item For all $x\in X$, $\langle x,\widetilde{\Phi}(\om)\rangle=\langle x,\Phi_0(\om)\rangle$, for $\mu$-a.e. \el
\end{lem}
\begin{proof}
Define $r(\om)=\|\Phi_0(\om)\|$ and hence $r\in L^q(\mu)$. Also, we have $\|[\Phi_0]\|_{q,w^*}=\|r\|_q$. Now, let $h\in L^q(\mu)$ be positive and finite everywhere with $h\geq r$. For a scalar lifting $\rho:\mc{L}^\iy(\mu)\ra\mc{L}^\iy(\mu)$ as stated above, define $$\widetilde{\Phi}(\om)(x)=h(\om)\rho (\Phi_x)(\om),$$ 
where $\Phi_x(\om)=\frac{\Phi_0(\om)(x)}{h(\om)}$. It is easy to check that $|\Phi_x|\leq \|x\|$ for all $x\in X$. This follows as a scalar function, $\Phi_x$ satisfies $0\leq |\Phi_x|\leq \chi_\Omega\|x\|$. Moreover, $$\rho(|\Phi_x|)(\om)\leq \rho(\chi_\Omega\|x\|)(\om)=\|x\|\rho(\chi_\Omega)(\om)=\|x\|.$$

{\sc Claim~:} There exists a countable set $\{x_1,...,x_n,...\}\subseteq X$ such that $$r(\om)=\sup_n\{|\Phi_0(\om)(x_n)|\},\quad \mu-a.e.$$

Let $\phi=I([\Phi_0])$ and choose $f_n\in B_{L^p(\mu,X)}$ such that $|\phi(f_n)|\to \|\phi\|$. Then the essential range of these functions $f_n$ lies in a common separable closed subspace $X_0\subseteq X$. Hence $\|\phi\|\geq \|\phi|_{L^p(\mu,X_0)}\|\geq \lim\limits_{n\to\iy}|\phi(f_n)|=\|\phi\|=\|r\|_q$. Choose $(x_k)$ dense in $B_{X_0}$ and set $$s(\om)=\sup_k\{|\Phi_0(\om)(x_k)|\}.$$
Then $s$ is measurable and $s(\om)=\|\Phi_0(\om)|_{X_0}\|\leq r(\om)$. Now for $f\in L^p(\mu,X_0)$, the definition of $s$ gives $$|\Phi_0(\om)(f(\om))|\leq s(\om)\|f(\om)\|\;\;\mu-a.e.$$ Now, we have 
    $$\begin{aligned}
        |\phi(f)|&=\Big|\int_\Om \Phi_0(\om)(f(\om)) d\mu(\om)\Big|\\
        &\leq \int_\Om |\Phi_0(\om)(f(\om))| d\mu(\om)\\
        &\leq \int_{\Om}s(\om)\|f(\om)\|d\mu(\om)\\
        &\leq \Big(\int_\Om s(\om)^q d\mu(\om)\Big)^\frac{1}{q}\Big( \int_\Om \|f(\om)\|^pd\mu(\om)\Big)^\frac{1}{p}\\
        &= \|s\|_q\|f\|_p.
    \end{aligned}$$
    Hence by taking the supremum over $f\in B_{L^p(\mu,X_0)}$, we get $$\|\phi\|=\|\phi|_{L^p(\mu,X_0)}\|\leq \|s\|_q\leq \|r\|_q=\|\phi\|.$$
    Since $0\leq s\leq r$ and $\|s\|_q=\|r\|_q$, the nonnegative function $r^q-s^q$ has integral zero, we have $s=r$ a.e. and consequently $r(\om)=\sup\limits_k\{|\Phi_0(\om)(x_k)|\}\quad \mu-a.e.$
    
Now, from the linearity of lifting and the definition of $\widetilde{\Phi}$, one can observe that $x\mapsto \widetilde{\Phi}(\om)(x)$ is linear for every $\om$. Moreover, $|\widetilde{\Phi}(\om)(x)|=h(\om)|\rho(\Phi_x)(\om)|\leq h(\om)\rho(|\Phi_x|)(\om)\leq h(\om)\|x\|$. Hence $\widetilde{\Phi}(\om)\in X^*$. Also, $\widetilde{\Phi}$ is $w^*$-measurable, which follows from the way it is defined. So for each $x\in X$, the lifting property (\cite[Definition 1.5.1.1]{CM}) gives $$\widetilde{\Phi}(\om)(x)=\Phi_0(\om)(x)\quad \mu-a.e.$$
    This follows that $[\widetilde{\Phi}]=[\Phi_0]$.

    Now, we establish the norm equality. For each $x\in B_X$, we have $$|\Phi_x|\leq \frac{r}{h}. $$
    So by lifting, we get $|\widetilde{\Phi}(\om)(x)|\leq h(\om)\rho(\frac{r}{h})(\om)$. Now taking the supremum over $B_X$ and using
$\rho(r/h)=r/h$ almost everywhere, yields $\|\widetilde{\Phi}(\om)\|\leq r(\om)$ $\mu$-a.e.

For the reverse inequality, first recall that $[\widetilde{\Phi}]=[\Phi_0]$ and the countable many vectors $\{x_1,...,x_k,..\}\subseteq X$ associated to $r(\om)$. Now, $$\begin{aligned}
    \|\widetilde{\Phi}(\om)\|&\geq \sup_k|\widetilde{\Phi}(\om)(x_k)|\\
    &=\sup_k|\Phi_0(\om)(x_k)|,\;\; \mu-a.e.\\
    &=r(\om),\;\; \mu-a.e.
\end{aligned}$$ 
Hence we get $\|\widetilde{\Phi}(\om)\|=\|\Phi_0(\om)\|$, $\mu$-a.e. and this completes the proof.
\end{proof}
The following lemma shows that the lifting construction in
Lemma~\ref{LL} preserves a prescribed strongly measurable
representative outside a single null set. This is stronger than
scalar equivalence and will be needed in Theorem~\ref{T12}.

\section{Stability of unique Hahn-Banach extension properties}
\begin{thm}\label{TU}
Let $(\Om,\Sigma,\mu)$ be a finite measure space and $Y$ be a subspace of $X$. Then $Y$ has property-$(U)$ in $X$ if and only if $L^p(\mu,Y)$ has property-$(U)$ in $L^p(\mu,X)$ for $1<p<\iy$.
\end{thm}
\begin{proof}
 Suppose that $Y$ has property-$(U)$ in $X$. Let $\psi(\neq 0)\in L^p(\mu,Y)^*$ and let $\phi_1,\phi_2\in L^p(\mu,X)^*$ be norm-preserving extensions of $\psi$. Assume that $\psi\neq 0$. As discussed earlier, write $$\psi=I([\Psi])\text{ and }\phi_i=I([\Phi_i]),\;i=1,2.$$
 For simplicity of notation again choose norm-minimal representatives by $\Psi,\Phi_1,\Phi_2,$ and put $$r_0(\om)=\|\Psi(\om)\|,\quad r_i(\om)=\|\Phi_i(\om)\|,i=1,2.$$
 Then,
 \begin{equation}\label{eqL0}
     \|r_0\|_q=\|[\Psi]\|_{q,w^*}=\|\psi\|=\|\phi_i\|=\|[\Phi_i]\|_{q,w^*}=\|r_i\|_q.
 \end{equation}

 Now fix $y\in Y$. Then for any $A\in\Si$, $$\begin{aligned}
     \int_A \Phi_i(\om)(y)d\mu(\om)&=\phi_i(y\chi_A)\\
     &= \psi (y \chi_A)\\
     &= \int_A \Psi(\om)(y)d\mu(\om)
 \end{aligned}$$
 This implies that $$\Phi_i(\om)(y)=\Psi(\om)(y),\quad \mu-a.e.$$
 Since $y$ is arbitrary, we get 
 \begin{equation}\label{eqL1}
     [\Phi_i|_Y]=[\Psi]\quad \text{in }L_{w^*}^q(\mu,Y^*).
 \end{equation}
 Now define $h=1+r_0+r_1+r_2$, then $h\in L^q(\mu)$ and positive everywhere. Now apply the Lemma~\ref{LL} using the same $h$ and same lifting to $\Psi,\Phi_1,\Phi_2$. Denote the corresponding elements by $\widetilde{\Psi}, \widetilde{\Phi_1}, \widetilde{\Phi_2}$. The above lemma gives 
 \begin{equation}\label{eqL2}
     [\widetilde{\Psi}]=[\Psi],\quad [\widetilde{\Phi_i}]=[\Phi_i]
 \end{equation} and 
 \begin{equation}\label{eqL3}
     \|\widetilde{\Psi}(\om)\|=r_0(\om),\quad \|\widetilde{\Phi_i}(\om)\|=r_i(\om),\quad \mu-a.e.
 \end{equation}
 For each $y\in Y$, equation(\ref{eqL1}) implies $\frac{\Phi_i(\om)(y)}{h(\om)}=\frac{\Psi(\om)(y)}{h(\om)}$, $\mu$-a.e. Then by lifting, we get \begin{equation}\label{eqL4}
     \widetilde{\Phi_i}(\om)|_Y=\widetilde{\Psi}(\om),\text{ for every }\om.
     \end{equation}
     Now, by (\ref{eqL4}), we have $\|\widetilde{\Psi}(\om)\|\leq \|\widetilde{\Phi_i}(\om)\|$. Using (\ref{eqL3}), we get $r_0(\om)\leq r_i(\om)$ $\mu$-a.e. Equation~(\ref{eqL0}) now gives $$\int_\Om \big(r_i(\om)^q-r_0(\om)^q\big)d\mu=\|r_i\|_q^q-\|r_0\|_q^q=0.$$
     This follows that $r_i=r_0$, $\mu$-a.e. Consequently, we get 
     \begin{equation}\label{eqL5}
        \|\widetilde{\Psi}(\om)\|=\|\widetilde{\Phi_1}(\om)\|=\|\widetilde{\Phi_2}(\om)\|,\; \mu-a.e. 
     \end{equation}
     Now, by (\ref{eqL4}) and (\ref{eqL5}), almost everywhere, $\widetilde{\Phi_1}(\om)$ and $\widetilde{\Phi_2}(\om)$ are norm-preserving extensions of $\widetilde{\Psi}(\om)$. Since $Y$ has property-$(U)$ in $X$, we get $$\widetilde{\Phi_1}(\om)=\widetilde{\Phi_2}(\om),\quad \mu-a.e.$$
     Therefore $[\widetilde{\Phi_1}]=[\widetilde{\Phi_2}]$ and by (\ref{eqL2}), we get $$[\Phi_1]=[\widetilde{\Phi_1}]=[\widetilde{\Phi_2}]=[\Phi_2].$$
 Finally, applying $I$, we obtain $\phi_1=\phi_2$, and it follows that $L^p(\mu,Y)$ has property-$(U)$ in $L^p(\mu,X)$.

     The converse follows from \cite[Remark~3.5(a)]{DPR}. This completes the proof.
\end{proof}
We now have the following corollary.
\begin{cor}\label{C6}
Let $(\Om,\Si,\mu)$ be a finite measure space, $Y$ be a closed subspace of a Banach space $X$ and $1<p<\iy$ with $\frac{1}{p}+\frac{1}{q}=1$. Suppose $\psi\in L^p(\mu,Y)^*$ and $\phi\in L^p(\mu,X)^*$ be a norm-preserving extension of $\psi$. Then there exist $w^*$-measurable norm-minimal representatives $\widetilde{\Psi}\in \mc{L}_{w^*}^q(\mu,Y^*)$ and $\widetilde{\Phi}\in \mc{L}_{w^*}^q(\mu,X^*)$ associated to $\psi$ and $\phi$ respectively such that $$\widetilde{\Phi}|_Y(\om)=\widetilde{\Psi}(\om)\text{ for every }\om\text{ and }\| \widetilde{\Phi}(\om)\|=\|\widetilde{\Psi}(\om)\|,\quad \mu-a.e.$$
\end{cor}

We now focus on the stability of property-$(SU)$. Recall that $Y$ has property-$(SU)$ in $X$ if and only if $Y$ has property-$(U)$ in $X$ and the Hahn-Banach extension map $T:Y^*\to X^*$ is linear. Here $T(y^*)=\{x^*\in X^*:\|x^*\|=\|y^*\|\text{ and }x^*|_Y=y^*\}$ for $y^*\in Y^*$. In this case, we also keep in mind the use of scalar lifting to obtain suitable representatives, as discussed in the proofs of Theorem~\ref{TU} and Lemma~\ref{LL}. 
\begin{thm}\label{TSU}
Let $(\Om,\Sigma,\mu)$ be a finite measure space and $Y$ be a subspace of $X$. Then $Y$ has property-$(SU)$ in $X$ if and only if $L^p(\mu,Y)$ has property-$(SU)$ in $L^p(\mu,X)$ for $1<p<\iy$.
\end{thm}
\begin{proof}
    Suppose $Y$ has property-$(SU)$ in $X$ and $T: Y^*\to X^*$ is the linear Hahn-Banach extension operator. Then from Theorem~\ref{TU}, we have $L^p(\mu,Y)$ has property-$(U)$ in $L^p(\mu,X)$. Now let $S: L^p(\mu,Y)^*\to L^p(\mu,X)^*$ be the Hahn-Banach extension operator. 
    
   {\sc Claim~:} $S$ is linear.

    Now, let $\psi_1,\psi_2\in L^p(\mu,Y)^*$ and $\psi_3= \psi_1+\psi_2$. Now put $\phi_i=S(\psi_i)$, $1\leq i\leq 3$. Clearly, $\phi_i$ are unique norm-preserving extensions of $\psi_i$, respectively.

     Now, choose the norm-minimal $w^*$-measurable representatives $$\Psi_i:\Om\to Y^*,\quad \Phi_i:\Om\to X^*$$ of $\psi_i$ and $\phi_i$ respectively. Here we also have for each $y\in Y$, $\Psi_3(\om)(y)=\Psi_1(\om)(y)+\Psi_2(\om)(y)$ $\mu$-a.e. as $\psi_3=\psi_1+\psi_2$. Define $h(\om)=1+\sum_{i=1}^3 \big(\|\Psi_i(\om)\|+\|\Phi_i(\om)\|\big)$, then $h\in L^q(\mu)$ and positive everywhere. Now apply Lemma~\ref{LL} using the same $h$ and the same lifting to all $\Psi_i$ and $\Phi_i$. Denote the corresponding elements by $\widetilde{\Psi_i},\widetilde{\Phi_i}$. Finally from the proof of Lemma~\ref{LL} and Theorem~\ref{TU}, we have the followings:
     \bla
      \item $\widetilde{\Phi_i}$ represents $\phi_i$ and $\widetilde{\Psi_i}$ represents $\psi_i$.
      \item $\widetilde{\Phi_i}(\om)|_Y=\widetilde{\Psi_i}(\om)$ for all $\om\in\Om$.
      \item $\widetilde{\Psi_3}(\om)=\widetilde{\Psi_1}(\om)+\widetilde{\Psi_2}(\om)$ for all $\om\in\Om$.
      \item $\|\widetilde{\Phi_i}(\om)\|=\|\widetilde{\Psi_i}(\om)\|$ $\mu$-a.e. $i=1,2,3$.
     \el
     From the above, it is clear that there exists $N\in\Si$ such that $\mu(N)=0$ and for $1\leq i\leq 3$, we have $\widetilde{\Phi_i}(\om)$ is a norm-preserving extension of $\widetilde{\Psi_i}(\om)$ for all $\om\in \Om\setminus N$. Now as $Y$ has property-$(U)$ in $X$, we have $\widetilde{\Phi_i}(\om)=T(\widetilde{\Psi_i}(\om))$. Observe that for $\om\in\Om\setminus N,$ $$\begin{aligned}
         \widetilde{\Phi_3}(\om)&=T(\widetilde{\Psi_3}(\om))\\
         &= T(\widetilde{\Psi_1}(\om)+\widetilde{\Psi_2}(\om))\\
         &= \widetilde{\Phi_1}(\om)+\widetilde{\Phi_2}(\om)\text{ [ Since $T$ is linear ]}
     \end{aligned}$$
      Finally, assume $f\in L^p(\mu,X)$, then $$\begin{aligned}
          S(\psi_1+\psi_2)(f)&=\phi_3(f)\\
          &= \int_\Om \widetilde{\Phi_3}(\om)(f(\om)) d\mu(\om)\\
          &= \int_\Om \widetilde{\Phi_1}(\om)(f(\om)) d\mu(\om)+\int_\Om \widetilde{\Phi_2}(\om)(f(\om)) d\mu(\om)\\
          &= S(\psi_1)(f)+ S(\psi_2)(f).
      \end{aligned}$$
      As $f\in L^p(\mu,X)$ is arbitrary, we get $S(\psi_1+\psi_2)= S(\psi_1)+S(\psi_2)$. Also it is easy to check that $S(\al \psi)=\al S(\psi)$ for all $\al\in \mb{R}$ and $\psi\in L^p(\mu,Y)^*$. Hence $S$ is linear, and this completes the proof.

      The converse follows from \cite[Remark~3.5(a)]{DPR}
\end{proof}
Now, we focus on the stability of property-$(wU)$. First, we recall the following theorem in this context from \cite[Theorem~5.2]{KPR}.
\begin{thm}
Let $(\Omega,\Sigma,\mu)$ be a probability space and $Y$ be a subspace of $X$. Suppose $X^*$ has the RNP, then $Y$ has property-$(wU)$ in $X$ if and only if $L^p(\mu,Y)$ has property-$(wU)$ in $L^p(\mu,X)$ for $1<p<\iy$.    
\end{thm}
In this article, the authors considered the question (\cite[Problem~5.3]{KPR}) of whether the RNP hypothesis can be weakened. We provide a positive answer without requiring the RNP. To begin, we note that the proof of the preceding theorem relies on \cite[Lemma~5]{DR}, which gives an equivalent characterization for a norm-attaining functional $\psi\in L^p(\mu,X^*)$, with $1<p<\iy$. That lemma, however, also uses the RNP assumption. Here we revisit these statements and establish them without imposing the RNP.
\begin{lem}\label{LNA}
    Let $(\Om,\Si,\mu)$ be a finite measure space, $X$ be a Banach space and $1<p<\iy$ with $\frac{1}{p}+\frac{1}{q}=1$. Also let $\phi\in L^p(\mu,X)^*$ and there exists $f\in S_{L^p(\mu,X)}$ such that $\phi(f)=\|\phi\|$. If $\Phi\in L_{w^*}^q(\mu,X^*)$ be a norm-minimal representative of $\phi$, then for almost all $\om\in\Om$, either $\Phi(\om)=0$ or $\Phi(\om)$ attains its norm.
\end{lem}
\begin{proof}
    Let us denote $D=\{\om\in\Om: \|f(\om)\|>0\},\; r(\om)=\|\Phi(\om)\|$ and $s(\om)=\|f(\om)\|$. Then recall that $\|r\|_q=\|\phi\|$ and $\|s\|_p=1$. Now, we have $$\begin{aligned}
        \phi(f)&=\int_D \Phi(\om)(f(\om))d\mu(\om)\\
        &\leq \int_D \|\Phi(\om)\|\|f(\om)\|d\mu(\om)\\
        &=\int_D r(\om)s(\om)d\mu(\om)\\
        &\leq \|r\chi_D \|_q \|s\|_p\\
        &\leq \|r\|_q\\
        &=\|\phi\|\\
        &= \phi(f).
    \end{aligned}$$ 
    This implies that $\Phi(\om)(f(\om))=r(\om)s(\om)$ $\mu$-a.e. on $D$ and that $\Phi(\om)=0$ $\mu$-a.e. on $\Om\setminus D$. It follows that there exists a null set $N\in\Si$ such that $\Phi(\om)=0$ for all $\om\in (\Om\setminus D)\setminus N$. Now we see that for all $\om\notin N$, $\Phi(\om)\neq 0$ implies $s(\om)> 0$.

    For this let $\om\notin N$ and $\Phi(\om)\neq 0$. If possible let $s(\om)=0$. But this implies that $\om\notin D$. Consequently $\om\in(\Om\setminus D)\setminus N $. This follows that $\Phi(\om)=0$, a contradiction. Recall that $s(\om)=\|f(\om)\|$ and  $\Phi(\om)(f(\om))=r(\om)s(\om)$ $\mu$-a.e. on $D$, we have 
     $$\Phi(\om)\Big(\frac{f(\om)}{\|f(\om)\|}\Big)=\frac{r(\om)s(\om)}{s(\om)}=r(\om)=\|\Phi(\om)\|.$$
     This completes the proof.
\end{proof}
Now, we see the converse of the above lemma.
\begin{lem}
    Let $(\Om,\Si,\mu)$ be a finite measure space, $X$ be a Banach space and $1<p<\iy$ with $\frac{1}{p}+\frac{1}{q}=1$. Also, let $\phi\in L^p(\mu,X)^*$ and that $\Phi\in L_{w^*}^q(\mu,X^*)$ be the norm-minimal representative of $\phi$. Suppose there exists $f\in L^p(\mu,X)$ such that $\Phi=0$ almost everywhere on $\{\om: f(\om)=0\}$ and for almost everywhere on $\{\om:f(\om)\neq 0\}$, $$\Phi(\om)\Big(\frac{f(\om)} {\|f(\om)\|}\Big)=\|\Phi(\om)\|.$$ Then there exists $g\in S_{L^p(\mu,X)}$ such that $\phi(g)=\|\phi\|$.
\end{lem}
\begin{proof}
  First let $r(\om)=\|\Phi(\om)\|$ and hence $r\in L^q(\mu)$. Moreover, we have, $\|\phi\|=\|r\|_q$ and this follows that $$\int_\Om r(\om)^q d\mu(\om)=\|\phi\|^q.$$
  Now, define $u:\Om\to X$ by $$u(\om)=\Bigg\{\begin{array}{ll}
     \frac{f(\om)}{\|f(\om)\|},  & f(\om)\neq 0 \\
       0, & f(\om)=0.
  \end{array}$$
  It is easy to check that $u$ is measurable as $f$ is measurable.

 {\sc Claim :} $\Phi(\om)(u(\om))=r(\om),\;\mu-a.e.$

By assumption, we have $\Phi(\om)(u(\om))=\|\Phi(\om)\|=r(\om)$ almost everywhere on $\{\om:f(\om)\neq 0\}$. Also, almost everywhere on  $\{\om:f(\om)=0\}$, we have $\Phi(\om)(u(\om))=0=\|\Phi(\om)\|$ ( Since $\Phi=0$ almost everywhere on $\{\om: f(\om)=0\}$ ). Hence, we have $\Phi(\om)(u(\om))=r(\om),\;\mu-a.e.$

Now let $P=\{\om:\|\Phi(\om)\|>0\}$. Since $\Phi=0$ almost everywhere on $\{\om:f(\om)=0\}$, one can easily check that $f(\om)\neq 0$ almost everywhere on $P$. Hence $\|u(\om)\|=1$ $\mu$-a.e. on $P$.
Now, define $g:\Om\to X$ by 
$$g(\om)= \Big(\frac{r(\om)}{\|\phi\|}\Big)^{q-1}u(\om).$$
It is easy to check that $g$ is measurable and $\|g(\om)\|=\Big(\frac{r(\om)}{\|\phi\|}\Big)^{q-1}$ $\mu$-a.e. on $P$. On $P^c$, $r(\omega)=0$, and therefore $g(\omega)=0$. Consequently,
$$
\|g(\omega)\|
=\left(\frac{r(\omega)}{\|\varphi\|}\right)^{q-1}
\quad\mu\text{-a.e. on }\Omega.
$$ Also, we have 
\[
      \|g\|_p^p=\int_\Om \Big(\frac{r(\om)}{\|\phi\|}\Big)^{p(q-1)}d\mu(\om)
       =\frac{1}{\|\phi\|^q}\int_\Om r(\om)^q d\mu(\om)
      =1.
  \]
  Hence, $g\in S_{L^p(\mu,X)}$. Again, by linearity of $\Phi(\om)$, we have $$\begin{aligned}
  \Phi(\om)(g(\om))&= \Big(\frac{r(\om)}{\|\phi\|}\Big)^{q-1} \Phi(\om)(u(\om))\\
  &= \Big(\frac{r(\om)}{\|\phi\|}\Big)^{q-1} r(\om)\quad \mu-a.e.\\
  &= \frac{r(\om)^q}{\|\phi\|^{q-1}}.
  \end{aligned}$$
  To see that $\|\phi\|=\phi(g)$, observe that 
  \[
      \phi(g)=\int_\Om \Phi(\om)(g(\om))d\mu(\om)
      =\frac{1}{\|\phi\|^{q-1}}\int_\Om r(\om)^q d\mu(\om)
      = \frac{\|\phi\|^q}{\|\phi\|^{q-1}}=\|\phi\|.
\]
  This completes the proof.
\end{proof}
Now, we have the main theorem concerning property-$(wU)$.
\begin{thm}\label{TwU}
    Let $(\Om,\Sigma,\mu)$ be a probability space and $Y$ be a subspace of $X$. Then $Y$ has property-$(wU)$ in $X$ if and only if $L^p(\mu,Y)$ has property-$(wU)$ in $L^p(\mu,X)$ for $1<p<\iy$.
\end{thm}
\begin{proof}
    Suppose $Y$ has property-$(wU)$ in $X$. Let $\psi\in L^p(\mu,Y)^*$ and attains its norm on $f\in S_{L^p(\mu,Y)}$. Let $\phi_1,\phi_2\in L^p(\mu,X)^*$ be two norm-preserving extensions of $\psi$. Clearly, $\phi_1,\phi_2$ also attain their norm on $f$. Let $\Psi,\Phi_1$ and $\Phi_2$ are the norm-minimal representations of $\psi,\phi_1$ and $\phi_2$ respectively. Applying the common lifting construction used in the proof of Theorem~\ref{TU}, with the same lifting and the same dominating function for $\Psi,\Phi_1,\Phi_2$, we obtain norm-minimal representatives $\widetilde\Psi,\widetilde\Phi_1,\widetilde\Phi_2$ satisfying the following properties. 
    \bla
      \item $\widetilde{\Phi_i}$ and $\widetilde{\Psi}$ are norm-minimal representatives of $\phi_i$ and $\psi$ respectively.
      \item $\widetilde{\Phi_i}(\om)|_Y=\widetilde{\Psi}(\om)$, for all $\om\in\Om$.
      \item $\|\widetilde{\Phi_i}(\om)\|=\|\widetilde{\Psi}(\om)\|$, $\mu$-a.e. $i=1,2$.
     \el
     Applying Lemma~\ref{LNA} to $\psi\in L^p(\mu,Y)^*$ and its norm-minimal representative $\widetilde\Psi$, we see that, almost everywhere non-zero $\widetilde\Psi(\omega)$ attains its norm on $Y$. For almost all $\om$ such that $\widetilde{\Psi}(\om)\neq 0$, properties (b) and (c), together with property-$(wU)$ of $Y$ in $X$, imply
$$
\widetilde\Phi_1(\omega)=\widetilde\Phi_2(\omega).
$$
If $\widetilde\Psi(\omega)=0$, condition (c) implies that both extensions are zero. Therefore $\widetilde\Phi_1=\widetilde\Phi_2$, $\mu$-a.e., and consequently $\phi_1=\phi_2$.

For the converse, we use the equivalent condition in Theorem~\ref{CwU}(c) of property-$(wU)$. For this let $x_1^*,x_2^*\in NA(X)$ with $x_i^*(y_i)=\|x_i^*\|$ and $x_1^*+x_2^*\in Y^\perp$, where $y_i\in S_{Y}$. We need to prove $x_1^*+x_2^*=0$. Now, observe that $\chi_\Om x_i^*\in L^q(\mu,X^*)\subseteq L^p(\mu,X)^*$ and $\chi_\Om x_i^*(\chi_\Om y_i)=\|\chi_\Om x_i^*\|$. Since $(x_1^*+x_2^*)(y)=0$ for all $y\in Y$, we have $\chi_\Om x_1^*+\chi_\Om x_2^*\in L^q(\mu,Y^\perp)\subseteq L^p(\mu,Y)^\perp$. As $L^p(\mu,Y)$ has the property-$(wU)$ in $L^p(\mu,X)$, we have $\chi_\Om x_1^*+\chi_\Om x_2^*=0$ and it follows that $x_1^*+x_2^*=0$. This completes the proof.
\end{proof}

\section{Examples on class representatives}

The following example shows that the pointwise extension depends on the choice of representatives. For a nonempty set $\Omega$, we denote the space $c_0(\Omega)$ to represent those $\zeta:\Omega\ra \mb{R}$ that satisfy $|\zeta (\om)|\geq \e$ for only finitely many $\omega$, for an $\e>0$.

\begin{ex}\label{ex1}
Let $\Om=[0,1]$ be the space endowed with the discrete topology and consider the Lebesgue measure space $([0,1],\Si,\mu)$. Define $X=\mb{R}\oplus_1 c_0(\Om)\oplus_1\mb{R}$ and $Y=\mb{R}\oplus_1 c_0(\Om)\oplus_1\{0\}$. Thus $\|(a,z,b)\|_X=|a|+\|z\|_\iy+|b|$. Also any $z\in c_0(\Om)$, one can observe that $z(\om)=0$ $\mu$-a.e.

Now for $\om\in \Om$, define $\Psi(\om)\in Y^*$ and $\Phi(\om)\in X^*$ by 
    
\begin{equation}\label{eq1}
     \Psi(\om)(a,z,0)=a\text{ and }\Phi(\om)(a,z,b)=a+z(\om).   
\end{equation}
    
Then one can check that $\|\Psi(\om)\|=\|\Phi(\om)\|=1$ and that the functions $\om\mapsto\Psi(\om)$ and $\om\mapsto\Phi(\om)$ are $w^*$-measurable. As $\mu$ is Lebesgue measure, one can also check that the norm functions $\om\mapsto\|\Psi(\om)\| \text{ and }\om\mapsto \|\Phi(\om)\|$ are in $L^q(\mu)$. Which in other words $\|\|\Psi(.)\|\|_q=\|\|\Phi(.)\|\|_q=1$.

Now, for $f\in L^p(\mu,X)$, write $f(\om)=(a_f(\om), z_f(\om),b_f(\om))$ and define $$\phi(f)=\int_\Om a_f(\om) d\mu(\om).$$
Now let $\psi=\phi|_{L^p(\mu,Y)}$. Then $\|\phi\|=\|\psi\|=1$. Now it is obvious that $\Psi$ represents $\psi$ and $\Phi$ represents $\phi$ in the sense of Theorem~\ref{T10}. Also, we have $\phi$ is a norm-preserving extension of $\psi$.

Now fix $\om\in \Om$ and $e_\om\in c_0(\Om)$ be defined by $$e_\om(s)=\Big\{\begin{array}{cc}
        1, & s=\om, \\
         0,& s\neq \om.
    \end{array}$$
Now, for $y_\om=(0,e_\om,0)$, equation~(\ref{eq1}), gives $$\Psi(\om)(y_\om)=0,\quad \Phi(\om)(y_\om)=1.$$
Therefore, $$\Phi(\om)|_Y\neq \Psi(\om)\text{  for each }\om\in\Om.$$
On the other hand, for fixed $y=(a,z,0)\in Y$, $\Phi(\om)(y)-\Psi(\om)(y)=z(\om)=0$, for $\mu$-a.e. Thus we have $$[\Phi|_Y]=[\Psi]\text{ but } \Phi(\om)|_Y\neq \Psi(\om)\text{ for all }\om\in\Om.$$
\end{ex}

\begin{ex}\label{ex2}
Let $\Om=[0,1]$ be the space endowed with the discrete topology and consider the Lebesgue measure space $([0,1],\Si,\mu)$.  Define $X= \mb{R}\oplus_1 c_0(\Om)$ and consider the space $L^p(\mu,X)$. Now, for $f\in L^p(\mu,X)$, write $f(\om)=(a_f(\om),z_f(\om))$. Define $\psi\in L^q(\mu,X^*)$ by $$\psi(f)=\int_\Om a_f(\om)d\mu(\om).$$ Then, $\psi$ has two norm-minimal representatives $\Psi$ and $\Phi$ with $[\Phi]=[\Psi]$ such that $\Psi$ is strongly measurable and $\Phi$ is $w^*$-measurable which is not strongly measurable.
\end{ex}
\begin{proof}
    For $\om\in \Om$, define $\Psi(\om)\in X^*$ and $\Phi(\om)\in X^*$ by 
    
    $$
     \Psi(\om)(a,z)=a\text{ and }\Phi(\om)(a,z)=a+z(\om).   
    $$
    From Example~\ref{ex1}, it is clear that both $\Psi,\Phi$ are norm-minimal representative of $\psi$ with $[\Psi]=[\Phi]$. Also, it is easy to check that $\Psi$ is Bochner measurable and $\Phi$ is $w^*$-measurable. Moreover, note that $$\|\Phi(\om_1)-\Phi(\om_2)\|=2\text{ for }\om_1\neq \om_2.$$
This implies that the range of $\Phi$ is not essentially norm-separable and hence $\Phi$ is not strongly measurable.
\end{proof}

\section{$w^*-w$-PC of unit norm functionals in $L^p(\mu,X)^*$}

In this section, we mainly study the sufficient condition on $\psi\in L^q(\mu,X^*)$ to be a $w^*-\sigma(L^p(\mu,X)^*,L^p(\mu,X^{**})$-PC of $B_{L^p(\mu,X)^*}$. However, we are unable to conclude the same for general $\psi\in L^p(\mu,X)^*$. Moreover, we see that this investigation is highly dependent on the choice of a representative of $\psi$. 
We now raise the following question, which we answer in an non-affirmative sense subsequently.

\begin{prob}
 Let $X$ be a Banach space with $X^*$ having RNP, $(\Omega, \Sigma, \mu)$ be a measure space and $1 < p < \infty$, $\frac{1}{p} + \frac{1}{q} = 1$. Let $\psi \in S_{L^q(\mu, X^*)}$ and $\Psi$ be a norm-minimal representative of $\psi$. If $\frac{\Psi(\om)}{\|\Psi(\om)\|}$ is a $w^*$-$w$-PC of $B_{X^*}$ $\mu$-a.e. with $\Psi(\om)\neq 0$, can one conclude that $\psi$ is a $w^*$-$\sigma(L^q(\mu, X^*), L^p(\mu, X^{**}))$-PC of $B_{L^q(\mu, X^*)}$?   
\end{prob}

Below is an example that serves as a counter example of the problem.

\begin{ex}\label{ex3}
Let $X, \psi,\Psi$ and $\Phi$ be as in Example~\ref{ex2}. Then every $\Phi(\om)$ is $w^*-w$-PC but no $\Psi(\om)$ is $w^*-w$-PC of $B_{X^*}$. Moreover $\psi$ is not a $w^*-w$-PC of $B_{L^q(\mu,X^*)}$.
\end{ex}
\begin{proof}
  First, recall that $\Psi(\om)=(1,0)$ and $\Phi(\om)=(1,e_\om)$, where $$ \Psi(\om)(a,z)=a\text{ and }\Phi(\om)(a,z)=a+z(\om)$$ and $e_\om(z)=z(\om)$ for all $z\in c_0(\Om)$.

  {\sc Claim 1:} Every $\Phi(\om)$ is $w^*-w$-PC of $B_{X^*}$.

  Let $\om\in \Om$ and $(a_\al,b_\al)\in B_{X^*}$ be a net such that $$(a_\al,b_\al)\stackrel{w^*}\longrightarrow \Phi(\om)=(1,e_\om).$$
  This follows that $a_\al\to 1$ and $b_\al(\om)\to 1$. Since, $\|b_\al\|_1\leq 1$, we have $$\begin{aligned}
      \|b_\al-e_\om\|_1&= |b_\al(\om)-1|+\sum_{\ga\neq \om}|b_\al(\ga)|\\
      &\leq |b_\al(\om)-1|+ 1-|b_\al(\om)|\\
      &\leq 2|b_\al(\om)-1|.
  \end{aligned}$$
  Hence, we have $b_\al\to e_\om$ and this follows that $(a_\al,b_\al)\stackrel{\|.\|}\longrightarrow (1,e_\om)$. Therefore $\Phi(\om)$ is a $w^*-\|.\|$-PC and this implies $\Phi(\om)$ is $w^*-w$-PC. Hence the claim follows.

  {\sc Claim 2:} No $\Psi(\om)$ is $w^*-w$-PC of $B_{X^*}$.

  Choose distinct $\om_n\in\Om$ and hence $z(\om_n)\to 0$ for all $z\in c_0(\Om)$. Consequently, we have $$(1,e_{\om_n})\stackrel{w^*}\longrightarrow (1,0)=\Psi(\om).$$
  Note that for $\chi_\Omega:\Omega\ra\mb{R}$, $x^{**}=(0,\chi_\Omega)\in X^{**}$ satisfies $$x^{**}(1,e_{\om_n})=1\text{ but }x^{**}(1,0)=0.$$
  This follows that $(1,e_{\om_n})\stackrel{w}\nrightarrow (1,0)$. Hence $\Psi(\om)$ is not a $w^*-w$-PC of $B_{X^*}$ and the claim is proved.

  {\sc Claim 3:} $\psi$ is not a $w^*-w$-PC of $B_{L^q(\mu,X^*)}$.

  Choose distinct $\om_n\in\Om$ and define $\psi_n\in B_{L^q(\mu,X^*)}$ by $$\psi_n(f)=\int_\Om \big(a(\om)+z(\om)(\om_n)\big)d\mu(\om), $$
  $\text{ where }f\in L^p(\mu,X)$ and $f(\om)=\big(a(\om),z(\om)\big)$. It is clear that $\psi_n$ is represented by the strongly measurable function $\psi_n$, defined by $\psi_n(\om)=(1,e_{\om_n})$. Now observe that $$|a(\om)+z(\om)(\om_n)|\leq |a(\om)|+\|z(\om)\|_\iy=\|f(\om)\|_X.$$
  Hence, from Hölder's inequality, we have $\|\psi_n\|\leq 1$. By choosing $f=(1,0)$, one can see $\|\psi_n\|=1$. Let $f\in L^p(\mu,X)$, where
$$
f(\om)=\bigl(a_f(\om),z_f(\om)\bigr),
$$
and define
$$
v_f=\int_\Om z_f(\om)d\mu(\om).
$$
Since the canonical projection from $\mb{R}\oplus_1 c_0(\Om)$ onto $c_0(\Om)$ is bounded and linear, the function $\om\mapsto z_f(\om)$ belongs to $L^p(\mu,c_0(\Om))$. Moreover, since $\mu(\Om)<\infty$, Hölder's inequality implies that $\om\mapsto z_f(\om)$ also belongs to $L^1(\mu,c_0(\Om))$. Therefore, its Bochner integral is well defined and belongs to $c_0(\Om)$. It follows that
$$
v_f=\int_\Om z_f(\om)d\mu(\om)\in c_0(\Om).
$$
Writing $v_f=(v_f(\ga))_{\ga\in\Om}$, fix $\ga\in\Om$ and consider the continuous coordinate functional
$$
\de_\ga:c_0(\Om)\longrightarrow\mb{R},\qquad
\de_\ga(z)=z(\ga).
$$
Since bounded linear operators commute with the Bochner integral, we obtain
\beqa
v_f(\ga)
=\de_\ga(v_f)
&=& \de_\ga\left(\int_\Om z_f(\om)d\mu(\om)\right)\\
&=& \int_\Om\de_\ga\bigl(z_f(\om)\bigr)d\mu(\om)\\
&=& \int_\Om z_f(\om)(\ga)d\mu(\om).
\eeqa
\[
\mbox{Now,}~~ \psi_n(f)-\psi(f)=\int_\Om z_f(\om)(\om_n)d\mu(\om)=v_f(\om_n)\ra 0,
\]
since $v_f\in c_0(\Om)$ and $\om_n$ are distinct. This follows that $\psi_n\stackrel{w^*}\longrightarrow \psi$.

  Now let $E:X\to L^p(\mu,X)$ be the canonical isometric embedding defined by $E(x)(\om)=x$. So $E^{**}:X^{**}\to L^p(\mu,X)^{**}$. Let $x^{**}=(0,I)\in X^{**}$ and hence $G=E^{**}x^{**}\in L^p(\mu,X)^{**}$, where $I(\om)=1$ for all $\om\in \Om$. By definition, $$G(\phi)=x^{**}(E^*\phi),\quad \phi\in L^q(\mu,X^*).$$
  Now, one can easily check that $E^*\psi_n=(1,e_{\om_n})$ and $E^*\psi=(1,0)$. This follows that $$G(\psi_n)=1\text{ for every }n\text{  and  }\quad G(\psi)=0.$$
  Hence $\psi$ is not a $w^*-w$-PC of $B_{L^q(\mu,X^*)}$ and the claim is proved.
\end{proof}
Now, we conclude about $w^*-w$-PC of elements in $L^q(\mu,X^*)$ of unit norms. First, we recall the following in our context.
In \cite{Rao}, the author introduces a stronger version of an ideal in Banach spaces, known as a {\it strict ideal}.
\bdfn 
A closed subspace $Y$ of a Banach space $X$ is called a {\it strict ideal} if there exists a norm one projection $P:X^*\ra X^*$
with  $\ker (P) = Y^{\perp}$ and  $B_{P(X^*)}$ is $w^*$-dense in $B_{X^*}$. 
\edfn
In {\cite[Lemma 1]{Rao}} it is established that a subspace $Y$ of $X$ is a strict ideal in $X$ if and only if $Y \ci X \ci Y^{**}$  under the canonical embedding of $Y$ in $Y^{**}$. An easy example of a strict ideal is a Banach space $X$ under its canonical embedding in $X^{**}.$ Following the natural inclusion $L^p(\mu,X)\subseteq L^p(\mu,X^{**})\subseteq L^p(\mu,X)^{**}$, it is now clear that $L^p(\mu,X)$ is a strict ideal in $L^p(\mu,X^{**})$, $1\leq p<\iy$.

\begin{rem}
Suppose that $X$ is a strict ideal of a Banach space $Z$, fix its linear Hahn-Banach extension operator $H:X^*\to Z^*.$ Then
it satisfies $$H(x^*)|_X=x^*,\quad \|H(x^*)\|=\|x^*\|\text{ and }\overline{H(B_{X^*})}^{\si(Z^*,Z)}=B_{Z^*}.$$

Then it is evident that 
$x_\al^*\stackrel{\si(X^*,Z)}\longrightarrow x^*\iff H(x_\al^*)\stackrel{\si(Z^*,Z)}\longrightarrow H(x^*)$. This fact is useful in our next observation.
\end{rem}

We need the following theorem from \cite[Proposition~2.1]{MP} which generalizes \cite[Lemma 2.14]{HWW}.
\begin{thm}\label{une}
Let $X$ and $Z$ be Banach spaces where $X$ is a strict ideal of $Z$. Then \tfae.
    \bla 
    \item $x^*\in S_{X^*}$ has a unique Hahn-Banach extension to a functional on $Z$.
    \item The function $I:(B_{X^*},w^*)\to (B_{X^*},\si(X^*,Z))$ is continuous at $x^*$.
    \el
\end{thm}

Now, we have the following theorem. The $w^*$-topology on $L^q(\mu,X^*)$ is simply $\si(L^q(\mu,X^*),L_p(\mu,X))$.

\begin{thm}\label{T12}
Let $X$ be a Banach space, $(\Omega, \Sigma, \mu)$ be a measure space and $1 < p < \infty$, $\frac{1}{p} + \frac{1}{q} = 1$. Let $\psi \in S_{L^q(\mu, X^*)}$. If $\frac{\psi(\om)}{\|\psi(\om)\|}$ is $w^*$-$w$-PC of $B_{X^*}$ $\mu$-a.e. whenever $\psi(\om)\neq 0$, then $\psi$ is a $w^*$-$\sigma(L^p(\mu, X)^*, L^p(\mu, X^{**}))$-PC of $B_{L^p(\mu, X)^*}$. 
\end{thm}

\begin{proof}
Put $E=L^p(\mu,X)$ and $Z=L^p(\mu,X^{**})$. Under the canonical embeddings,
$$
E\subseteq Z\subseteq E^{**},
$$
so $E$ is a strict ideal in $Z$. Therefore, by Theorem~\ref{une}, it is enough to prove that $\psi$, regarded as an element of $E^{*}$, has a unique norm-preserving extension to $Z$.
Define $\widehat{\psi}\in Z^{*}$ by
$$
\widehat{\psi}(u)=\int_{\Om}u(\om)\bigl(\psi(\om)\bigr)d\mu(\om),
\qquad u\in Z.
$$
By Holder's inequality, $\|\widehat{\psi}\|\leq\|\psi\|_q=1$. Moreover, for every $f\in E$,
$$
\widehat{\psi}(f)=\int_{\Om}\psi(\om)\bigl(f(\om)\bigr)d\mu(\om)
=\psi(f).
$$
Thus $\widehat{\psi}$ extends $\psi$, and consequently
$$
1=\|\psi\|\leq\|\widehat{\psi}\|\leq 1.
$$
Hence $\widehat{\psi}$ is a norm-preserving extension of $\psi$ to $Z$.
Since $\psi\in L^q(\mu,X^{*})$ is strongly measurable, its norm-minimal representative may be chosen to be $\psi$ itself. As $L^p(\mu,X)\ci L^p(\mu,X^{**})$, $\psi$ has a Hahn-Banach extension $\phi$ over $L^p(\mu,X^{**})$. 

Applying Corollary~\ref{C6} for the function $\Phi$, we obtain a norm-minimal $w^{*}$-measurable representative
$$
\Phi:\Om\longrightarrow X^{***}
$$
of $\phi$ such that
$$
\Phi(\om)|_X=\psi(\om)
$$
and
$$
\|\Phi(\om)\|=\|\psi(\om)\|
$$
for $\mu$-a.e.
Fix $\om$ outside the corresponding null set. If $\psi(\om)=0$, then the norm equality gives $\Phi(\om)=0$. Suppose that $\psi(\om)\neq0$. By hypothesis,
$$
\frac{\psi(\om)}{\|\psi(\om)\|}
$$
is a $w^{*}$-$w$ point of continuity of $B_{X^{*}}$. Since $X$ is a strict ideal in $X^{**}$, Theorem~\ref{une} shows that $\psi(\om)$ has a unique norm-preserving extension to $X^{**}$. Both $\Phi(\om)$ and the canonical image $J_{X^{*}}\psi(\om)$ are such extensions. Therefore,
$$
\Phi(\om)=J_{X^{*}}\psi(\om)
$$
for $\mu$-almost every $\om\in\Om$.
Consequently, for every $u\in Z$,
$$
\begin{aligned}
\phi(u)
&=\int_{\Om}\Phi(\om)\bigl(u(\om)\bigr)d\mu(\om)\\
&=\int_{\Om}J_{X^{*}}\psi(\om)\bigl(u(\om)\bigr)d\mu(\om)\\
&=\int_{\Om}u(\om)\bigl(\psi(\om)\bigr)d\mu(\om)\\
&=\widehat{\psi}(u).
\end{aligned}
$$
Thus $\phi=\widehat{\psi}$, and hence $\psi$ has a unique norm-preserving extension to $Z$. Another application of Theorem~\ref{une} now shows that $\psi$ is a $w^{*}\text{-}\sigma\bigl(L^p(\mu,X)^{*},L^p(\mu,X^{**})\bigr)
$ point of continuity of $B_{L^p(\mu,X)^{*}}$.
\end{proof}
However, we do not have a similar result for the converse of Theorem~\ref{T12}. 

\begin{prob}
Let $X$ be a Banach space where $X^*$ has RNP, $(\Omega, \Sigma, \mu)$ be a measure space and $1 < p < \infty$, $\frac{1}{p} + \frac{1}{q} = 1$. Let $\psi \in S_{L^q(\mu, X^*)}$ is a $w^*-w$-PC of $B_{L^q(\mu, X^*)}$ . Can we say $\frac{\psi(\om)}{\|\psi(\om)\|}$ is $w^*$-$w$-PC of $B_{X^*}$ $\mu$-a.e. with $\psi(\om)\neq 0$?     
\end{prob}
\begin{thm}\label{Lpstab}
Let $X$ be a Banach space, $(\Omega, \Sigma, \mu)$ a measure space and $1 < p < \infty$, $\frac{1}{p} + \frac{1}{q} = 1$. Let $(E_i)\ci \Si$ be a countable or finite family of pairwise disjoint sets with $\mu(A_i)>0$. Then $g=\sum\limits_{i=1}^{\infty} x_i^* \chi_{E_i} \in S_{L^p(\mu,X)^*}$ is a $w^*$-$\sigma(L^p(\mu,X)^*, L^p(\mu, X^{**}))$-PC of $B_{L^p(\mu,X)^*}$ if and only if $\frac{x_i^*}{\|x_i^*\|}$ is a $w^*$-$w$-PC of $B_{X^*}$ for all $i$.
\end{thm}
\begin{proof}
Suppose that $g=\sum\limits_{i=1}^{\infty} x_i^* \chi_{E_i} \in S_{L^p(\mu,X)^*}$ is a $w^*$-$\sigma(L^p(\mu,X)^*, L^p(\mu, X^{**}))$-PC of $B_{L^p(\mu,X)^*}$.
Let $(y_\alpha^*)$ be a net in $B_{X^*}$ such that  $y_\alpha^* \xrightarrow{w^*} x_k^* / \|x_k^*\|$ for some $k$, then define $g_\alpha = \sum_{i \neq k} x_i^* \chi_{E_i} + \|x_k^*\| y_\alpha^* \chi_{E_k} \in B_{L^p(\mu,X)^*}.$ Then $g_\alpha \xrightarrow{w^*} g$. 

By assumption, $g_\alpha \xrightarrow{\sigma(L^p(\mu,X)^*, L^p(\mu, X^{**}))} g$. 

Thus $\|x_k^*\| y_\alpha^* \chi_{E_k} \xrightarrow{\sigma(L^p(\mu,X)^*, L^p(\mu, X^{**}))} x_k^* \chi_{E_k}.$ 

So, $ y_\alpha^* \chi_{E_k} \xrightarrow{\sigma(L^p(\mu,X)^*, L^p(\mu, X^{**}))} \frac{x_k^*}{\|x_k^*\|} \chi_{E_k}.$ Hence, $y_\alpha^* \xrightarrow{w} x_k^* / \|x_k^*\|.$
The converse follows from Theorem~\ref{T12}.
\end{proof}
\begin{cor}
Let $X$ be a Banach space, $(\Omega, \Sigma, \mu)$ a measure space and $1 < p < \infty$, $\frac{1}{p} + \frac{1}{q} = 1$. Let 
$$\mathcal{A}=\{g\in S_{L^q(\mu, X^*)} : g=\sum_{i=1}^{\infty} x_i^* \chi_{E_i}, \ x_i^*\in X^*, E_i\in \Sigma, \forall i\}.$$ 
Then the following are equivalent.
    
\bla 
\item Each element of $\mathcal{A}$ has a unique Hahn-Banach extension to a functional on $L^p(\mu, X^{**}).$
\item $X$ is Hahn-Banach smooth.
\item $L^p(\mu,X)$ is Hahn-Banach smooth.
\item Each element of $\mathcal{A}$ has a unique Hahn-Banach extension to a functional on $L^p(\mu, X)^{**}.$
\el 
\end{cor}

\begin{proof}
$(a)\Ra(b).$  Let $x^*\in S_{X^*}.$ Then $x^*\chi_{\Omega}\in \mathcal{A}.$ Since $x^*\chi_{\Omega}$ has unique Hahn-Banach extension to a functional on $L^p(\mu, X^{**})$, what follows from Theorem~\ref{une} and the discussion above, 
   $$
   I: \Big(B_{L^p(\mu,X)^*}, \si\big(L^p(\mu,X)^*, L^p(\mu,X)\big)\Big)\to \Big(B_{L^p(\mu,X)^*},\; \si\big(L^p(\mu,X)^*, L^p(\mu,X^{**})\big)\Big) 
   $$
is continuous at $x^*\chi_{\Omega}$. 
Then by Theorem  $\ref{Lpstab}$, $x^*$ is $w^*-w$-PC. This follows that each $x^*\in S_{X^*}$ is $w^*-w$-PC and hence $X$ is Hahn-Banach smooth. 
   
$(b)\Ra (c)$ follows from Corollary~8 of \cite{HL}. 
   
$(c)\Ra (d)\Ra (a)$ is obvious.
\end{proof}

\end{document}